\documentclass[11pt, openany]{book}

\usepackage{amsmath}
\usepackage{amssymb}
\usepackage{latexsym}
\usepackage{amsxtra}
\usepackage{amscd}
\usepackage{theorem}
\usepackage{sectsty}
\usepackage{calrsfs}

\chapterfont{\Large}
\sectionfont{\large}

\theoremstyle{plain}

{\theorembodyfont{\slshape}

	\newtheorem{thm}{Theorem}[chapter]
	\newtheorem{theorem}[thm]{Theorem}
	\newtheorem{question}[thm]{Question}
	\newtheorem{prob}[thm]{Problem}
	\newtheorem{problem}[thm]{Problem}
	
	\newtheorem{proposition}[thm]{Proposition}
	\newtheorem{corollary}[thm]{Corollary}

	\newtheorem{fact}[thm]{Fact}
	\newtheorem{qu}[thm]{Question}

}

{\theorembodyfont{\rmfamily}

	\newtheorem{remark}[thm]{Remark}
	\newtheorem{defi}[thm]{Definition}
	\newtheorem{example}[thm]{Example}

}

\newcommand{\n}{\noindent}
\newcommand{\vn}{\vspace{3mm} \noindent}
\newcommand{\Pic}{\mbox{\rm Pic}}
\newcommand{\dm}{\mbox{\rm dim}}

\newcommand{\B}{{\hfill$\Box$}}

\def\CC{{\mathbb C}}
\newcommand \CA {{\cal A}}
\newcommand \ZZ {{\mathbb Z}}
\newcommand \QQ {{\mathbb Q}}
\newcommand \CM {{\cal M}}
\newcommand \CH {{\cal H}}
\newcommand \GG {{\mathbb G}}
\newcommand \Spec {\mathop{\rm Spec}}
\newcommand \Jac {\mathop{\rm Jac}}
\newcommand  \FF {{\mathbb F}}
\newcommand \dime {\mathop{\rm dim}}
\newcommand \Hom {\mathop{\rm Hom}}
\newcommand \PP {{\mathbb P}^1}
\newcommand \Aa {{\mathbb A}^1}
\newcommand{\Z}{{\mathbb Z}}
\newcommand{\el}{\ell}
\newcommand{\Aut}{\mbox{\rm Aut}}
\newcommand{\Gal}{\mbox{\rm Gal}}
\newcommand{\PPP}{{\mathbb P}}
\newcommand{\F}{{\mathbb F}}
\newcommand{\N}{{\mathbb N}}
\newcommand{\Fr}{\mbox{\rm Fr}}
\newcommand{\HH}{\mathcal{H}}
\newcommand{\M}{\mathcal{M}}
\newcommand{\T}{\mathcal{T}}
\newcommand{\la}{\longrightarrow}

\def\cT{{\cal T}}

\def\cN{{\cal N}}
\def\cA{{\cal A}}

\def\cH{{\cal H}}
\def\cM{{\cal M}}

\begin{document}

\begin{titlepage}

\hrule 

\bigskip

{\Large Problems from the workshop on  \\

{\bf Automorphisms of Curves} \\ 

(Leiden, August, 2004)}

\vspace{1cm}

{\large {\sl edited by} Gunther Cornelissen and Frans Oort

\bigskip

{\sl with contributions of} I.~Bouw, T.~ Chinburg, G.~Cornelissen, C.~Gasbarri, D.~Glass, C.~Lehr, M.~Matignon, F.~Oort, R.~Pries,  S.~Wewers}

\bigskip

\hrule

\vspace{2cm}

\noindent In the week of August, 16th -- 20th of 2004, we organized a workshop about
``Automorphisms of Curves'' at the Lorentz Center in Leiden. The programme included two ``problem sessions''. Some of the problems presented at the workshop were written down; this is our edition of these refereed and revised papers. 

The editing process was simplified in that the bibliographies of consecutive papers were put together at the end, without combining them. Thus, some
references might occur several times on the (non-alphabetical) list.

We thank all contributors and (anonymous) referees for their fast reaction, which allowed us to present this timely text only two months after the end of the workshop. 

We also thank the participants of the workshop for creating a lively and creative atmosphere. 

We hope some of the problems proposed will stimulate research and will soon be solved!

\bigskip

Utrecht, November, 1st, 2004

\bigskip

Gunther Cornelissen

Frans Oort

\end{titlepage}

\setcounter{tocdepth}{0}
\tableofcontents

\chapter{Rational functions with given monodromy on generic curves \\ {\mdseries \textsl{by} Irene I.\ Bouw \textsl{and} Stefan Wewers}}

The suggested problems are concerned with finding rational functions
on curves with prescribed monodromy group and ramification behavior.
Typically, this kind of problem is difficult to solve for a fixed
curve, so one first tries to solve it for general curves.

In our recent paper \cite{An}, we extend a result of V\"olklein and
Magaard \cite{VoelkleinMagaard03} to positive characteristic. This
result says that the generic curve of genus $g$ admits a rational
function of degree $n$ with alternating monodromy group, for all but
finitely many values of $n$. Except in characteristic $2$ and $3$, we
are able to determine the minimal degree of such a function, depending
on $g$.

In characteristic $2$, the smallest case that we were not able to
handle is the following:

\begin{prob} \label{BWP1}
  Let $X$ be the generic curve of genus $1$ in characteristic
  $2$. Show that there exists a rational function
  $f:X\to \PPP^1$ of degree $5$ with monodromy group $A_5$ and with
  ramification type $(3,3,3,3,3)$. 
\end{prob}

We believe that such a function $f$ exists, because we have a good
candidate for it. Let $k_0$ denote an algebraic closure of $\FF_2$,
and let $E$ be the supersingular elliptic curve over $k_0$. Taking the
quotient of $E$ under an automorphism of order $3$, we obtain a
rational function $g:E\to\PPP^1_{k_0}$ which reveals $E$ as the cyclic
cover of degree $3$ of $\PPP^1$ with three branch points. Using
patching, one can construct a tamely ramified cover $f:X\to\PPP^1_k$ of
degree $5$ with $X$ of genus one, defined over an algebraically closed field $k$ of
transcendence degree $2$ over $k_0$, with the following properties (see
\cite{An}, Proposition 4):
\begin{itemize}
\item
  the monodromy of $f$ is $A_5$,
\item
  the ramification type of $f$ is $(3,3,3,3,3)$,
\item the branch points of $f$ are $0,1,\infty,\lambda,\mu\in\PPP^1_k$,
  where $\lambda,\mu$ are elements of $k$ which are algebraically
  independent over $k_0$,
\item
  for `$\lambda=\mu$', the cover $f$ degenerates to the cover
  $g:E\to\PPP_{k_0}^1$. 
\end{itemize}
The last point deserves further explanation. By this we mean that
there exists a valuation ring $R\subset k$ and an $R$-model
$f_R:X_R\to Z_R$ of $f$ such that the following hold: a.)
$\bar{\lambda}=\bar{\mu}\neq 0,1,\infty$ in the residue field
$k_1:=R/\mathfrak{m}$ of $R$, b.) $X_R$ and $Z_R$ are semistable
curves and $f_R$ is an admissible cover (see \cite{HarrisMorrison}, \S
3.G), c.) the special fiber $X_R\otimes k_1$ contains a unique
component $X_1$ of genus $1$, and d.) the restriction of $f_R$ to
$X_1$ can be identified with the cover $g\otimes k_1:E\otimes
k_1\to\PPP^1_{k_1}$.

The cover $f:X\to\PPP_k^1$ corresponds to a $5$-tuple
$\sigma=(\sigma_1,\ldots,\sigma_5)$ of $3$-cycles which generate $A_5$
and verify the relation $\prod_i\sigma_i=1$. It degenerates to the
cover corresponding to the $4$-tuple
$\sigma'=(\sigma_1,\sigma_2,\sigma_3,\sigma_4\sigma_5)$. However, by
construction we have $\sigma_4\sigma_5=1$, and hence we get a cyclic
cover of degree $3$ with three branch points on the special fiber.

To solve Problem \ref{BWP1} one only needs to show that the curve $X$ is
generic, i.e.\ not isotrivial. However, we have not been able to prove
this. One problem is that $X$ has good reduction at the valuation
corresponding to $R$. Therefore the main argument used in
\cite{VoelkleinMagaard03} and \cite{An} to show that a curve is
generic fails.

Note that the analogous statement (i.e.\ genericity of $X$) in
characteristic $0$ (and even in characteristic $p>5$) holds. This is
proved in \cite{FriedKlaKop00}, using topological arguments.

Here is a more general version of Problem \ref{BWP1}:

\begin{prob} \label{BWP2}
  Let $G\subset S_n$ be a transitive permutation group on $n\geq 3$
  letters and $\sigma=(\sigma_1,\ldots,\sigma_r)$ a tuple of
  generators of $G$, satisfying the relation $\prod_i\sigma_i=1$. Let
  $\HH(\sigma)$ denote the Hurwitz space over $\ZZ$ which
  parameterizes tamely ramified covers $f:X\to\PPP^1$ of degree $n$ and
  with `branch cycle description' $\sigma$ (see e.g.\ 
  \textup{\cite{FriedVoelklein91}}). Let
  \[
       \phi_\sigma: \HH(\sigma) \;\to\; \M_g
  \]
  be the map which sends the class of a cover $f:X\to\PPP^1$ to the
  class of the curve $X$. Let $p$ be a prime number ($p=0$ is also
  allowed).
  \begin{enumerate}
  \item
    Is $\HH(\sigma)\otimes_{\ZZ}\FF_p$ nonempty?
  \item
    If the answer to (1) is `yes', compute the dimension of the image of
    $\phi_\sigma\otimes\FF_p$.
  \end{enumerate}
\end{prob}

The classical case is $G=S_n$, with $\sigma$ an $r$-tuple of
transpositions and $r=2n+2g-2$. In this case, it is known that
$\HH(\sigma)\otimes_{\ZZ}\FF_p$ is nonempty if and only if $p\neq 2$.
Furthermore, if $p>n$ and $2n-2\leq g$ (resp.\ $2n-2\geq g$) then the
image of $\phi_\sigma\otimes\FF_p$ has dimension $2g+2n-5$ (resp.\ is
dense in $\M_g$). See \cite{ACGH} and \cite{Fulton69}.

Note that Problem \ref{BWP1} corresponds to the special case of Problem \ref{BWP2} with
$G:=A_5$, and where $\sigma$ is any $5$-tuple of $3$-cycles in $A_5$
with product one. The construction that we have sketched above shows
that $\HH(\sigma)\otimes_{\ZZ}\FF_p$ is nonempty for $p\not=3$. Also, in
characteristic $p>5$, one can show that $\phi_\sigma\otimes\FF_p$ has
a dense image. However, the argument breaks down for $p\leq 5$,
because the Hurwitz space $\HH(\sigma)$ has bad reduction at $2$, $3$
and $5$.

The results of \cite{VoelkleinMagaard03} and \cite{An} solve Problem \ref{BWP2}
in many cases, with $G=A_n$. All these results are obtained by `going
to the boundary' of $\HH(\sigma)$ and $\M_g$. Problem \ref{BWP1} is an example
where this method seems to fail. We suggest another method that one
could try:

\begin{prob} \label{BWP3}
  Let $G\subset S_n$ and $\sigma=(\sigma_1,\ldots,\sigma_r)$ be as in
  Problem 2. Let $f:X\to\PPP^1$ be a cover of type $\sigma$ which
  corresponds to a generic point of $\HH(\sigma)\otimes_{\ZZ}\FF_p$. Compute
  the rank of the map induced by $\phi_\sigma$ on the tangent spaces at the
  point corresponding to $f$:
  \[
      {\rm d}\phi_\sigma|_f:\T_{\HH(\sigma),f}\;\to\;\T_{\M_g,X}.
  \]
\end{prob}

Note that Problem \ref{BWP3} is essentially equivalent to Problem \ref{BWP2}.2.
However, the different formulation suggests a new approach. There is a
similar set of problems called {\em infinitesimal Torelli problems for
  generalized Prym varieties}. They come up when we associate to an
\'etale Galois cover $f:X\to Y$ the abelian subvariety of $J_X$
corresponding to an irreducible $\QQ$-representation of the Galois
group $G$. Here there are some recent results by Tamagawa
\cite{Tamagawa} which might be a guideline for Problem \ref{BWP3}.

\bigskip

{\small Irene Bouw, Heinrich Heine Universit\"at, Mathematisches Institut, Universit\"atsstra{\ss}e 1,  40225 D\"usseldorf, Deutschland, bouw@math.uni-duesseldorf.de

\medskip 

Stefan Wewers, Mathematisches Institut der
Rheinischen Friedrich-Wilhelms-Uni\-versit\"at Bonn, Beringstra{\ss}e 1, 53115 Bonn, Deutschland, wewers@math.uni-bonn.de}

\chapter{Can deformation rings of group representations not be local complete intersections? \\
{\mdseries \textsl{by} Ted Chinburg}}

Suppose $G$ is a profinite group and that $k$ is a field of positive characteristic $p$.  Let $V$
be a finite dimensional vector space over $k$ with the discrete
topology having a continuous $k$-linear action of $G$.  If
$\mathrm{End}_k(V) = k$, it is known by work of Schlessinger,
Mazur, Faltings, de Smit and Lenstra that $V$ has a universal
deformation ring $R(G,V)$.  (See \cite{dSL} for the defining
properties and a construction of $R(G,V)$;  we take the auxiliary ring 
$\mathcal{O}$ 
in \cite{dSL} to be the ring of infinite Witt vectors
over $k$.)  

The following question is due to Matthias Flach.

\begin{question}
{\rm } Are there $G$ and $V$ as above for which
$R(G,V)$ is Noetherian but
not a local complete intersection?
\end{question} 

Note that by the argument given in \cite[Thm. 2.3.3]{dSL} and at the end of the proof of 
Lemma 7.3 of \cite{dSL}, 
$R(G,V)$ will be Noetherian if and only if ${\rm dim}_k H^1(G,{\rm End}_k(V)) < \infty$.

The ring $R(G,V)$ has been shown in many cases to be a local complete
intersection.  As one example, suppose $k$ is finite, ${\rm dim}_k(V) = 2$
and that $G$ is the absolute Galois group $G_K$ of a local field 
$K$ of residue characteristic $p$.  In \cite{B}, G. B\"ockle shows
that $R(G,V)$ is a complete intersection which is flat over $W(k)$.

\bigskip

{\small Department of Mathematics, University of Pennsylvania, Philadelphia, USA, 

ted@math.upenn.edu}

\chapter{Lifting an automorphism group to finite characteristic \\ {\mdseries \textsl{by} Gunther Cornelissen}}

Let $X$ be a smooth projective curve over a field $k$ of characteristic $p>0$ and $\rho: G \hookrightarrow \mbox{Aut}(X)$ a finite group of automorphisms of $X$. Let $W(k)$ denote the Witt vectors of $k$. Let $D(X,\rho)$ denote the deformation functor of $(X,\rho)$ that assigns to an element $A$ of the category $\mathcal{C}$ of local artinian $W(k)$-algebras with residue field $k$ the set of isomorphism classes of liftings of $(X,\rho)$ to $A$. Here, a lifting of $(X,\rho)$ to $A$ is a smooth scheme $X'$ of finite type over $A$, an isomorphism of its special fiber with $X$, and an action $G \hookrightarrow \mbox{Aut}_A(X')$ that lifts $\rho$ in the obvious sense. Then $D(X,\rho)$ has a prorepresentable hull $R(X,\rho)$ in the sense
of Schlessinger (and is actually prorepresentable if $H^0(X,T_X)^G=0$, e.g., if the genus $g$ of $X$ satisfies $g \geq 2$, \cite{BM}, 2.1). 

Recall that the characteristic $\mbox{char}(R)$ of a ring $R$ is the positive generator of 
the kernel of the unique morphism ${\bf Z} \rightarrow R$. 
Let $\nu(X,\rho)$ denote the characteristic of $R(X,\rho)$. This number is a power of $p$ or zero.

\begin{problem}
Give an example of a group action $(X,\rho)$ such that $\nu(X,\rho) = p^n$ for $n>1$ finite. 
\end{problem}

\begin{remark} \label{ref1} If the second ramification group for the local action of $G$ at a point of $X$ is trivial, we say $(X,\rho)$ is weakly ramified locally at that point. We say that $(X,\rho)$ is weakly ramified if it is weakly ramified everywhere locally. For example, if $X$ is an ordinary curve, any group action on it is weakly ramified (S.\ Nakajima, \cite{Nak}). A calculation by myself and Ariane M\'ezard (\cite{CM}) gives the following result: \emph{if $(X,\rho)$ is weakly ramified and $\nu(X,\rho) \neq p$, then $\nu(X,\rho)=0$ (hence no example as above can be produced) and all non-trivial ramification groups of order divisible by $p$ are on the following list:} $${\bf Z}/p; D_p \mbox{ or $\mathrm{[}A_4$ and $p=2\mathrm{]}$.}$$ As an easy exercise, prove Oort's conjecture for weakly ramified
actions of a cyclic group. The above calculation tells you more: it classifies all weakly ramified group actions into liftable and non-liftable ones and gives
the corresponding hulls and the group action on them explicitly (\cite{BM} 4.2 for $\nu(X,\rho)=0$ and \cite{CK}, 4.1 for $\nu(X,\rho)=p$). One sees in particular that if $\nu(X,\rho)=0$, then $R(X,\rho)$ is flat over $W(k)$. 
\end{remark}

\begin{remark}
The calculation in remark (\ref{ref1}) appears to be compatible 
with forthcoming work of T.~Chinburg, D.~Harbater and R.~Guralnick that {\sl exclude}
certain group actions from being liftable to characteristic zero. They determine the set $\mathcal L$ of abstract groups $G$ with a normal $p$-Sylow subgroup $P$ and $G/P$ cyclic, such that there exists an embedding $G \hookrightarrow \mbox{Aut}_k (k[[x]])$, for which the local lifting obstruction of Bertin (\cite{Bertin}) is non-trivial (not necessarily weakly ramified). They ask whether the complement of this set $\mathcal L$ is exactly the set of groups $G$ for which {\sl every} embedding $G \hookrightarrow \mbox{Aut}_k (k[[x]])$ can be lifted. 

The set of groups in the first remark is contained in the complement of this set $\mathcal L$. The calculation with M\'ezard actually shows that {\sl no} locally weakly ramified action of a $G$ outside the list can lift beyond characteristic $p$, and {\sl any} locally weakly ramified action of $G$ on the list can be lifted to characteristic zero. Note that in general, liftability can depend on the action (not just the abstract type of the group): one can give an example of two embeddings of the same abstract group $G$ into $\mbox{Aut}_k (k[[x]])$, one of which lifts to characteristic zero, and one of which doesn't. Indeed,  
let $G=({\bf Z}/p)^2$. There is a (non-weakly ramified) action of $G$ on $k[[x]]$ that lifts to characteristic zero by the results of Matignon in \cite{Mat}. On the other hand, the action of $G$ on $k[[x]]$ by $ x \mapsto {x}/{(1+ux)}$ for $u$ running through a 2-dimensional ${\bf F}_p$-vector space in $k$ does not lift to characteristic zero.    
\end{remark}

\begin{remark}
Any action of a group on the above list (\ref{ref1}) lifts to characteristic zero (irrespective of being weakly ramified), as follows from combining the following results: Oort-Sekiguchi-Suwa (\cite{OSS}) deal with the case of ${\bf Z}/p$; Pagot (\cite{Pagot}) treats $D_2={\bf Z}/2 \times {\bf Z}/2$, and Bouw and Wewers treat $D_p$ (\cite{BW}) and $A_4$ for $p=2$.  
\end{remark}

\begin{remark}
M.~Matignon remarks that one might look at liftability obstructions coming from ``equidistant geometry conditions'' as in the work of G.\ Pagot. Such 
conditions might fail to hold at a finite (non-prime) characteristic.
\end{remark}

\begin{remark}
Assume $R(X,\rho)$ pro-represents $D(X,\rho)$. We have $$\nu(X,\rho) = \sup \{ \mbox{char}(A) \ : \ A \in L(X,
\rho) \}$$  with  $$L(X,\rho) = \{ A \in \mathcal{C}  : \mbox{ there exists a lift of $(X,\rho)$ to $A$} \},$$ where we agree that $\sup p^n$ for $n$ strictly increasing equals $0$. Indeed, every artinian quotient of $R(X,\rho)$ belongs to $L(X,\rho)$, so $\nu(X,\rho)$ is less than or equal the indicated supremum. On the other hand, if $(X,\rho)$ lifts to $A$ of characteristic $p^m$, then there is a morphism $R(X,\rho) \rightarrow A$ and hence $p^m \leq \nu(X,\rho)$. 
\end{remark}

\begin{remark}
Suppose furthermore that for $(X,\rho)$, $\nu(X,\rho) \in \{0,p\}$. Then if one can find {\sl one} element $A \in \mathcal{C}$ of characteristic $p^2$ to which $(X,\rho)$ lifts, then it automatically lifts to some ring of characteristic zero. [Note, however, that the category $\mathcal{C}$ can contain rings $A$ of characteristic $p^m$ that do not lift to rings $B$ of characteristic $p^{m+1}$ (in the sense that $B/p^m=A$). We can therefore not conclude from our hypothesis $\nu(X,\rho) \in \{0,p\}$ that any lift of $(X,\rho)$ to a ring of characteristic $p^2$ lifts
to a ring of characteristic $0$ in this sense.] The result in (\ref{ref1}) shows that for weakly ramified actions, $p^2$ is indeed the critical characteristic. If one could a priori prove that $\nu(X,\rho) \in \{0,p\}$ holds for any $(X,\rho)$, this would give another approach to liftability questions (such as Oort's conjecture); but there is no further reason to believe that $p^2$ is the critical value. 
\end{remark}
 
\bigskip

{\small Universiteit Utrecht, Mathematisch Instituut, Postbus 80.010, 3508 TA Utrecht, Nederland, cornelis@math.uu.nl}

\chapter{Flat connections and representations of the fundamental group in characteristic $p>0$ \\ {\mdseries \textsl{by} Carlo Gasbarri}}

Let $k$ be the algebraic closure of  a finite field. Let $X$ be a
connected smooth projective curve over $k$. Let $F\colon X\to X$
be the absolute Frobenius of $X$.

\bigskip

It is well known that there is an equivalence of categories
between:
\begin{itemize}
\item The category of vector bundles equipped with a stratification
$(E, \nabla)$ and horizontal morphisms. We recall that a
stratification is a morphism of Lie algebras $\nabla\colon
\textrm{Diff}(X)\to \mathrm{End}_k(E)$ such that, if $s$ is a local section of $E$,
$D$ a local section of $\mathrm{Diff}(X)$ and $f$ a local regular function
on $X$ then
$$\nabla(D)(fs)=f\nabla(D)(s)+\nabla(D_f)(s)$$
where $D_f$ is the differential operator $D_f(g):=D(fg)-fD(g)$
(operator of degree at most one less then the degree of $D$).

\item The category $SS(X)$ of sequences $\{ E_n , \sigma_n\}$ where
$E_n$ is a vector bundle over $X$ and $\sigma_n\colon
F^\ast(E_n)\buildrel\sim\over\longrightarrow E_{n-1}$ is an
isomorphism.
\end{itemize}

We can prove that:

\begin{fact} Given a representation $\rho: \pi_1(X)\to GL_N(k)$ of the
fundamental group of $X$, we can associate to it a vector bundle
equipped with a stratification $E_\rho$; moreover the sequence
corresponding to it is such that the set $\{ (E_n, \sigma_n)\}$ is \emph{finite}. 
\end{fact} 

In order to prove this it suffices to remark that
the image of $\rho$ is contained in $GL_N(\F_q)$ for some $q$
(because $\pi_1(X)$ is finitely generated) so it is finite: this
implies that, if $F:X\to X$ is the absolute Frobenius then
$F^n(E_\rho)\simeq E_\rho$.

\begin{fact}
Given a sequence $\{ E_n, \sigma_n\}$ as above such that the
set $\{ (E_n, \sigma_n)\}$ is finite, the corresponding stratified
vector bundle comes from a representation of the fundamental group
of $X$. 
\end{fact}

This essentially follows from \cite{LS}.

\bigskip

Now, for every positive integer $r$, we can consider the category
$SS(X)_r$, of which the objects are $\{ E_n, \sigma_n\}$ with
$\mathrm{Card}(\{ (E_n, \sigma_n)\})\leq r$.
We can prove that:

\begin{fact} $SS(X)_r$ is a Tannakian category.
\end{fact} 

The tensor
product structure is the tensor product between vector bundles
with a stratification and morphisms are horizontal morphisms (cf.\ \cite{Gi}), the functor fibre is the functor ``restriction to a closed
point''; so it is the category of the representations of a group
$\pi_1^{(r)}(X)$. In order to prove that $SS(X)_r$ is Tannakian, we
essentially work as in \cite{No}: it is easy to see that the vector
bundles $E_n$ are semistable of degree zero so the functor fibre
is fully faithful etc.\

\bigskip

By abstract nonsense there is a surjection
$\pi_1(X)\to\pi_1^{(r)}(X)$; moreover, if $r>r'$ we have a
surjection $\pi_1^{(r)}(X)\to \pi_1^{(r')}(X)$.

We think that it is true that
$$\pi_1(X)=\lim_\leftarrow\pi_1^{(r)}(X).$$

\begin{question} How much of the $\pi_1^{(r)}(X)$'s determine $X$?\end{question}

\begin{question} Can we characterize the $\pi_1^{(r)}(X)$'s geometrically:  that
means without any mention of stratified vector bundles?\end{question}

\noindent By this we mean that we can characterize the Galois
coverings of $X$ whose group surjection $\pi_1(X)\to G$
factors through $\pi_1^{(r)}(X)$.

\begin{fact} For every $N$, each $\pi_1^{(r)}(X)$ has only
finitely many irreducible representations of rank $N$.
\end{fact}

\begin{fact} If $X$ is an elliptic curve we have that
$$\pi_i^{(r)}(X)\simeq T_p(X)\times X[p^r-1]$$ \end{fact}

This can be computed by using the classification of stratified
vector bundles on an elliptic curve given in \cite{Gi}, page 10.

\begin{question} Are the $\pi_1^{(r)}(X)$'s extensions of a finite group by a
$p$--group? \end{question}

Or, inspired by the genus one case:

\begin{question}
How is the finiteness of $\pi_1^{(r)}(X)$ related with the fact that the Jacobian of $X$ is (or is not) ordinary?
\end{question}

And, perhaps more important, a (meta)question: \emph{are the $\pi_1^{(r)}(X)$'s interesting? Are they useful for
something?}

\bigskip

{\small Dipartimento di Matematica dell'Universit\`a di
Roma II ``Tor Vergata", Viale della Ricerca Scientifica, 00133
Roma (Italy), gasbarri@mat.uniroma2.it}

\chapter{Questions on $p$-torsion of hyperelliptic curves \\ {\mdseries \textsl{by} Darren Glass \textsl{and} Rachel Pries}}

\section{Introduction}
We describe geometric questions raised by recent work on 
the $p$-torsion of Jacobians of curves defined over an algebraically
closed field $k$ of characteristic $p$.
These questions involve invariants of the $p$-torsion such as the $p$-rank or $a$-number.  
Such invariants are well-understood and have been
used to define stratifications of the moduli space $\CA_g$ of
principally polarized abelian varieties of dimension $g$.  
A major open problem is to understand how the Torelli locus intersects such strata in $\CA_g$. 
In \cite{GP}, we show that some of these strata intersect the image of the hyperelliptic locus under the Torelli map.  
This work relies upon geometric results on the configurations of branch points for 
non-ordinary hyperelliptic curves and raises some new geometric questions.

\section{Notation}
Let $k$ be an algebraically closed field of characteristic $p$.
Consider the moduli space $\CM_g$ (resp.\ $\CH_g$) of smooth (resp.\ hyperelliptic) curves of genus $g$.

The group scheme $\mu_p=\mu_{p,k}$ is the kernel of Frobenius on $\GG_m$, so
$\mu_p \simeq \Spec(k[x]/(x-1)^p)$.
If $\Jac(X)$ is the Jacobian of a $k$-curve $X$, the {\it $p$-rank}, $\dime_{\FF_p} \Hom(\mu_p, \Jac(X))$, of $X$
is an integer between $0$ and $g$.
A curve of genus $g$ is said to be {\it ordinary} if it has $p$-rank equal to $g$.
In other words, $X$ is ordinary if $\Jac(X)[p] \cong (\ZZ/p \oplus \mu_p)^g$.
Let $V_{g,f}$ denote the sublocus of curves of genus $g$ with $p$-rank at most $f$.
For every $g$ and every $0 \leq f \leq g$, the locus $V_{g,f}$ has codimension $g-f$ in
$\overline{\CM}_g$, \cite{FVdG}.  

The group scheme $\alpha_p=\alpha_{p,k}$ is the kernel of Frobenius on $\GG_a$,
so $\alpha_p \simeq \Spec(k[x]/x^p)$.
The {\it $a$-number}, $\dime_k \Hom(\alpha_p, \Jac(X))$, of $X$ is an integer 
between $0$ and $g$.
A generic curve has $a$-number equal to zero. A supersingular
elliptic curve $E$ has $a$-number equal to one.  In this case
there is a non-split exact sequence $0 \to \alpha_p \to E[p] \to
\alpha_p \to 0$.  
There is a unique isomorphism type of group scheme for 
the $p$-torsion of a supersingular elliptic curve, which we denote $M$. 
Let $T_{g,a}$ denote the sublocus of curves of genus $g$ with $a$-number at least $a$.

Let $N$ be the group scheme corresponding to the $p$-torsion of a 
supersingular abelian surface which is not superspecial.
By \cite[Example A.3.15]{G:book}, there is a filtration
$H_1 \subset H_2 \subset N$ where $H_1 \simeq \alpha_p$,
$H_2/H_1 \simeq \alpha_p \oplus \alpha_p$ and $N/H_2 \simeq \alpha_p$.
Moreover, the kernel $G_1$ of Frobenius and the kernel $G_2$ of Verschiebung are contained in $H_2$
and there is an exact sequence $0 \to H_1 \to G_1 \oplus G_2 \to H_2 \to 0$.
Finally, let $Q$ be the group scheme corresponding to the $p$-torsion of an
abelian variety of dimension three with $p$-rank $0$ and $a$-number $1$.   

These group schemes can be described in terms of their covariant Dieudonn\'e modules.
Consider the non-commutative ring $E=W(k)[F,V]$ with the Frobenius automorphism $\sigma:W(k) \to W(k)$
and the relations $FV=VF=p$ and $F\lambda=\lambda^\sigma F$ and $\lambda V=V \lambda^\sigma$ for all
$\lambda \in W(k)$.
Recall that there is an equivalence of categories between finite commutative group schemes
${\mathbb G}$ over $k$ (with order $p^r$) and finite left $E$-modules $D({\mathbb G})$
(having length $r$ as a $W(k)$-module).
By \cite[Example A.5.1-5.4]{G:book}, $D(\mu_p)=k[F,V]/k(V, 1-F)$, $D(\alpha_p)=k[F,V]/k(F,V)$,
and $D(N)=k[F,V]/k(F^3,V^3,F^2-V^2)$.
One can also show that $D(Q)=k[F,V]/k(F^4,V^4,F^3-V^3)$.

The $p$-rank of a curve $X$ with $\Jac(X)[p] \simeq N$ is zero.
To see this, note that $\Hom(\mu_p, N)=0$ or that
$F$ and $V$ are both nilpotent on $D(N)$.
The $a$-number of a curve $X$ with $\Jac(X)[p] \simeq N$ is one
since $N[F] \cap N[V] =H_1 \simeq \alpha_p$.

\section{Results}

Here are the results from \cite{GP} on the $p$-torsion of hyperelliptic curves.

\begin{theorem}\label{Thyperprk}
For all $0 \leq f \leq g$, the locus $V_{g,f} \cap \CH_g$ is non-empty of dimension $g-1+f$.
In particular, there exists a smooth hyperelliptic curve of genus $g$ and $p$-rank $f$.
\end{theorem}

The proof follows from the fact that $V_{g,0} \cap \CH_g$ is non-empty \cite{FVdG},
the purity result of \cite{OdJ}, and a dimension count at the boundary of $\CH_g$.

For the rest of the paper, suppose $p > 2$.
We consider the sublocus $\CH_{g,n}$ of the moduli space $\CM_g$ consisting
of smooth curves of genus $g$ which admit an action by $(\ZZ/2\ZZ)^n$ so that the quotient is the projective line.
We analyze the curves in the locus $\CH_{g,n}$ 
in terms of fibre products of hyperelliptic curves. 
We extend results of Kani and Rosen
\cite{KR} to compare the $p$-torsion of the Jacobian of a curve $X$ in
$\CH_{g,n}$ to the $p$-torsion of the Jacobians of its $\ZZ/2\ZZ$-quotients
(up to isomorphism rather than up to isogeny).

This approach allows us to produce families of Jacobians of (non-hyperellip\-tic) curves
whose $p$-torsion contains interesting group schemes.  The difficulty lies in 
controling the $p$-torsion of all of the hyperelliptic quotients of $X$.
This reduces the study of $\Jac(X)[p]$ to the 
study of the intersection of some subvarieties in
the configuration space of branch points.
For example, for Corollaries \ref{Tanumb2} and \ref{Canum3}, we study the 
geometry of the subvariety defined by Yui corresponding to the 
branch loci of non-ordinary hyperelliptic curves, \cite{Y}.  
Similarly, we use this method to show 
that $T_{g,a} \cap \CM_g$ is non-empty under certain conditions on $g$ and $a$.

In special cases, these families of curves intersect $\CH_g$.
This leads to the following partial results on the existence of hyperelliptic curves 
with interesting types of $p$-torsion.

\begin{corollary} \label{Cdim2}
Let $N$ be the $p$-torsion of a supersingular abelian surface which is not superspecial.
For all $g \geq 2$, there exists a smooth hyperelliptic curve $X$ of genus $g$ so that
$\Jac(X)[p]$ contains $N$.
\end{corollary}

Corollary \ref{Cdim2} is proved inductively starting with a curve $X$ of genus $2$ with $\Jac(X)[p] =N$.
In fact, we expect the $p$-torsion of the generic point of $V_{g,g-2} \cap \CH_g$ (which has dimension $2g -3$)
to have group scheme $N \oplus (\ZZ/p \oplus \mu_p)^{g-2}$.
We explain in \cite{GP} how this would follow from an affirmative answer to Question \ref{Q1}. 
 
\begin{corollary} \label{Tanumb2} Suppose $g \geq 2$ and $p \geq 5$.
There exists a dimension $g-2$ family of smooth
hyperelliptic curves of genus $g$ whose fibres have $p$-torsion 
$M^2 \oplus (\ZZ/p \oplus \mu_p)^{g-2}$
(and thus have $a$-number equal to $2$).
\end{corollary}

In fact, we expect $T_{g,2} \cap V_{g,g-2} \cap \CH_g$ to have dimension $2g-4$.  

\begin{corollary} \label{CQ}  
Let $Q$ be the $p$-torsion of an abelian variety of dimension three with $p$-rank $0$ and $a$-number $1$.
Suppose $g \geq 3$ is not a power of two.  
Then there exists a smooth hyperelliptic curve $X$ of genus $g$ so that
$\Jac(X)[p]$ contains $Q$.
\end{corollary}

Corollary \ref{CQ} is proved inductively starting from the supersingular hyperelliptic 
curve $X$ of genus 3 and $a$-number $1$ (and thus $\Jac(X)[p]=Q$) from \cite{O}.
An affirmative answer for $g=4$ in Question \ref{Q2} would allow us to remove the restriction on $g$ in 
Corollary \ref{CQ}.  We expect the generic point of $V_{g,g-3} \cap \CH_g$ (which has dimension $2g -4$)
to have group scheme $Q \oplus (\ZZ/p \oplus \mu_p)^{g-3}$.

\begin{corollary} \label{Canum3} Suppose $g \geq 5$ is odd and $p \geq 7$.
There exists a dimension $(g-5)/2$ family of smooth
hyperelliptic curves of genus $g$ whose fibres have $p$-torsion containing $M^3$
(and thus $a$-number at least $3$).
\end{corollary}

To determine the precise form of the group scheme in Corollary \ref{Canum3}, one could
consider Question \ref{Qlast}.
We expect $T_{g,3} \cap V_{g,g-3} \cap \CH_g$ to have dimension $2g-7$.  

\section{Questions}

These results raise the following geometric questions.
First, the expectations in the preceeding section all rest on the assumption that 
natural loci such as $\CH_g$, $V_{g,f}$ and $T_{g,a}$ should intersect as transversally as possible.
This tranversality can be measured both in terms of dimension and tangency.  
The meaning behind Theorem \ref{Thyperprk} is that the intersection of $\CH_g$ and $V_{g,f}$ at least 
has the appropriate dimension.  We could further ask this question.  

\begin{question} Is $V_{g,f} \cap \CH_g$ reduced? 
\end{question}

This relates to the question of whether $\CH_g$ and $V_{g,f}$ are transversal in the strict 
geometric sense.  Our proof of Corollaries \ref{Tanumb2} and \ref{Canum3} required us to show that
$V_{g,g-1} \cap \CH_g$ is not completely non-reduced.

An affirmative answer to the next question would imply by our fibre product construction
that for all $g \geq 4$ there exists a smooth hyperelliptic curve $X$ with
$\Jac(X)[p] \simeq N \oplus (\ZZ/p \oplus \mu_p)^{g-2}$.  This would then be the $p$-torsion
of the generic point of $V_{g,g-2} \cap \CH_g$.  

\begin{question} \label{Q1}
Given an arbitrary hyperelliptic cover $C \to \PP_k$, is it possible to
deform $C$ to an ordinary hyperelliptic curve by moving only one of the branch points?
\end{question}

The answer to Question \ref{Q1} will be affirmative if 
the hypersurface studied by \cite{Y} does not contain a line parallel to a coordinate axis.
To rephase this, consider the branch locus $\Lambda=\{\lambda_1,\ldots, \lambda_{2g}\}$
of an arbitrary hyperelliptic curve of genus $g-1$.  Does there
exist $\mu \in \Aa_k -\Lambda$ so that the hyperelliptic curve branched at
$\{\lambda_1,\ldots,\lambda_{2g}, \infty, \mu \}$ is ordinary?
For a generic choice of $\Lambda$, the answer to this question is yes, but this is 
not helpful for interesting applications.

One would like to strengthen Corollary \ref{Canum3} to state that there are 
hyperelliptic curves of genus $g$ with $a$-number exactly three.  
This raises a question which we state here in its simplest case.
Recall that $\lambda$ is supersingular if the elliptic curve branched at 
$\{0,1,\infty, \lambda\}$ is supersingular.  There are $(p-1)/2$ supersingular values of $\lambda$ by \cite{I}.

\begin{question} \label{Qlast} 
Which of the group schemes $(\ZZ/p \oplus \mu_p)^2$, $M \oplus (\ZZ/p \oplus \mu_p)$,
$N$, or $M^2$ occur as the $p$-torsion of the hyperelliptic curve branched at 
$\{0,1,\infty, \lambda_1, \lambda_2, \lambda_3\}$ when $\lambda_1, \lambda_2, \lambda_3$ are distinct supersingular values?
\end{question}

We expect that for all $p$ there exist distinct supersingular values $\lambda_1, \lambda_2, \lambda_3$,
so that the hyperelliptic curve branched at $\{0,1,\infty, \lambda_1, \lambda_2, \lambda_3\}$ is ordinary.  
This has been verified by C. Ritzenthaler for $7 \leq p < 100$.  
An affirmitive answer for any $p$ implies that there exists a smooth hyperelliptic curve of genus $5$ in characteristic $p$, 
with $p$-rank $2$ and $a$-number $3$.  
For multiple values of $p$, it appears that there do not exist distinct supersingular values $\lambda_1, \lambda_2, \lambda_3$,
so that the hyperelliptic curve branched at $\{0,1,\infty, \lambda_1, \lambda_2, \lambda_3\}$ has $p$-rank $0$. 

\begin{question} \label{Q2}
Does there exist a smooth hyperelliptic curve $X$ of genus $4$ (resp. $5$) so that 
$\Jac(X)[p]=Q \oplus (\ZZ/p \oplus \mu_p)$ (resp.\ $Q \oplus (\ZZ/p \oplus \mu_p)^2$)?
\end{question}

One would guess that the answer to Question \ref{Q2} is yes, but there does not seem to be much data.
If the answers to Questions \ref{Q1} and \ref{Q2} are both affirmative, then the generic point 
of $V_{g,g-3} \cap \CH_g$ has $p$-torsion of the form $Q \oplus (\ZZ/p \oplus \mu_p)^{g-3}$ for all $g \geq 3$.

\bigskip

{\footnotesize The first author was partially supported by NSF VIGRE grant DMS-98-10750.
The second author was partially supported by NSF grant DMS-04-00461.}

{\small

\bigskip

Darren Glass,
Department of Mathematics,
Columbia University,
New York, NY 10027, USA,
glass@math.columbia.edu

\medskip

Rachel J. Pries,
101 Weber Building,
Colorado State University,
Fort Collins, CO 80523-1874, USA,
pries@math.colostate.edu
}

\chapter{Automorphisms of curves and stable reduction \\ {\mdseries \textsl{by}  Claus Lehr \textsl{and} Michel Matignon}}

\section{$p$-groups as automorphism groups of curves in characteristic $p$}

The reference is \cite{Le-Ma1} and \cite{Le-Ma2}, 
the second of which is still in progress.
By 
$k$ we denote an algebraically closed field of characteristic $p>0$.

\subsection{Rationality conditions}

Let $C_f: W^p-W=f(X)\in k[X]$ be a $p$-cyclic cover of the affine 
line over $k$.
In \cite{Le-Ma1}, we have shown that a $p$-Sylow subgroup of the automorphism 
group of $C_f$ is an 
extensions of $\Z/p\Z$ by an elementary abelian $p$-group $(\Z/p\Z)^n$
for $n\geq 0$. Here the kernel of the extension, $\Z/p\Z$, is the group
of the Artin-Schreier cover.
Conversely any such group extension occurs as $p$-Sylow of the
automorphism group of some $C_f$ in the above maner.

\begin{qu} 
For a given group extension as above is it 
possible to find $f(X)\in \F_p[X]$ such that $C_f$ realizes this 
extension? (We want the
kernel of the extension to be the group of the Artin-Schreier
cover given by $C_f$).
\end{qu}

\noindent {For example if the group $D_8$ is the $2$-Sylow
of the automorphism group of the curve $C_f$ then one shows that 
$\deg f=1+\el 2^s$ for some 
$\el>1$ prime-to-$p=2$ and $s\geq 3$. Hence the minimal value of this degree
is $1+ 3.2^3=25$.
We gave such an $f\in \F_{16}[X]$ (\cite{Le-Ma1} 7.2 B. ii)) and an example 
with $f\in \F_2[X]$  (of degree $1+5.2^3=41$) is given in section 5.3 there.} 

\subsection{Big $p$-group actions: Nakajima condition (N)}

\begin{defi}  Let  $C$ be a non singular projective curve over $k$ of genus 
$g_C$ and  $G$ a $p$-subgroup of $\Aut_k C$.
We say that $(C,G)$ satisfies condition \textrm{\textbf{(N)}} if $$\textbf{(N)} \ \ \ \ g_C>0 \mbox{ and }\frac{|G|}{g_C}>\frac{2p}{p-1}.$$
\end{defi}

This definition is motivated by the following proposition from 
\cite{Le-Ma1} which is a translation of 
results  in \cite{Na}.

\begin{proposition} Assume $(C,G)$ with $g_C\geq 2$ satisfies condition \emph{\textbf{(N)}}.
Then there is a point, say  $\infty\in C$, such that $G$ is the wild inertia subgroup of $G$ at $\infty$.
Moreover
$C/G$ is isomorphic to $\PPP^1_k$ and the ramification locus (resp. branch locus) of the cover $\pi: C\to C/G$ is the point $\infty$ (resp. $\pi(\infty)$). We denote the ramification groups in lower 
numbering by $G_i$ ($i\geq 0$).
Let $i_0$ be the integer such that $G_2=G_3=....=G_{i_0}\supsetneqq G_{i_0+1}$. Then

i) we have $G_2\neq G_1$ and the quotient curve $C/G_2$ is isomorphic to $\PPP^1_k$.

ii) if we let $H$ be a subgroup which is normal in $G$ and such that $g_{C/H}>0$, then $G/H$ is a $p$-subgroup of $\Aut_k C/H$ and $$ \frac{|G|}{g_C}\leq \frac{|G/H|}{g_{C/H}}.$$ In particular $(C/H,G/H)$ satisfies 
condition \emph{\textbf{(N)}}. Moreover if $M\leq \frac{|G|}{g_C^2}$ for some $M$ one gets 
$$|H|\leq \frac{1}{M}\frac{|G/H|}{g_{C/H}^2}.$$

iii) if we let $H$ be a subgroup which is normal in $G$ and $G_2\supsetneqq H\supset G_{i_0+1}$, then $$g_{C/H}=\frac{(|G_2/H|-1)(i_0-1)}{2}>0$$ and $(C/H,G/H)$ satisfies condition \emph{\textbf{(N)}}.  

\end{proposition}

\begin{qu} 
Can one classify the pairs $(C,G)$ with 
$g_C\geq 2$ and $G$ a $p$-group of automorphisms of $C$ satisfying
condition \emph{\textbf{(N)}}? 
\end{qu}

\noindent
The following are some more precise questions. {\bf We shall always assume
$(C,G)$ satisfies condition (N).}

\begin{qu} 
What are the possible groups that can occur for $G_2$?   
\end{qu}

\noindent{It seems there are serious restrictions on $G_2$. 
In all examples we know 
$G_2$ is abelian and of exponent $p$. So we can ask}

\begin{qu}
Is it possible for $G_2$ to be non abelian?
\end{qu}

\noindent The following is a serious restriction on $G_2$.

\begin{theorem}[\cite{Le-Ma2}] Let $(C,G)$ satisfy \emph{\textbf{(N)}}. Then $G_2=G'$ is the 
commutator subgroup (it is also equal to $G'G^p$, the Frattini subgroup of 
$G$). 
\end{theorem}

\begin{corollary}
If $G_2$ is non-abelian then $Z(G_2)$ cannot be cyclic.
\end{corollary}

\noindent For example a non-abelian group of order $p^3$ and, more generally, extraspecial $p$-groups cannot occur for $G_2$.

\begin{corollary}
Let $(C,G)$ satisfy \emph{\textbf{(N)}} and assume that $|G_2|=p^3$. Then 
$G_2$ is abelian. 
\end{corollary}

\begin{qu}
Assume $G_2$ is abelian. Is it possible for the exponent of $G_2$ to be greater than $p$?
\end{qu}

\noindent {We can prove that $G_2$ cyclic implies $G_2\simeq \Z/p\Z$ 
(cf. \cite{Le-Ma2}).
One can ask if there is an action 
with $G_2=(\Z/p^2\Z)\times (\Z/p\Z)^t$
and, if the answer is yes, to give a bound on $t$.}

\begin{qu}\
Classify the actions where $G_2$ is elementary abelian.
\end{qu}
In \cite{Le-Ma2} we have solved the case where $G_2\subset Z(G)$ and we have given 
examples with $G_2\not\subset Z(G)$.

\begin{qu}
Describe the set of possible values of $\frac{|G|}{g_C^2}$ and for each 
value give a bound 
on the dimension of the subvariety in $M_g$ corresponding to curves $C$ with
such a big action. (Recall that we consider actions satisfying 
condition (N)).
\end{qu}

\noindent 
{In \cite{Le-Ma1} we proved that if $ \frac{|G|}{g_C^2}\geq
 \frac{4}{(p-1)^2}$ then only two values for  $\frac{|G|}{g_C^2}$
are possible, namely: $$\frac{4}{(p-1)^2} \mbox{ and }\frac{4p}{(p-1)^2}.$$
The last case corresponds to the actions $(C_f, G_{\infty,1})$ where $f=XR(X)$ with $R(X)$ an additive polynomial.
The group 
$G_{\infty,1}$ is then an extraspecial group and the dimension  of the subvariety in $M_g$ corresponding to curves $C$ with such a big action is O($\log(g)$) 
where $g=\frac{(p-1)}{2}p^{s}$. In \cite{Le-Ma2} we give such a description 
for  $$ \frac{|G|}{g_C^2}\geq \frac{4p^3}{(p^3-1)^2}.$$}

\section{Semi-stable reduction of $p$-group covers of 
curves over $p$-adic fields}

Intimate relations between characteristics $0$ and $p>0$ are encoded in
the geometry of Galois covers of curves 
over a discrete valuation ring $R$ of unequal characteristic $(0,p)$. 
One challenge is the desingularisation of such covers and to describe the 
action of the Galois group $\Gal (K^{alg}/K)$ on semi-stable models. 

The problem becomes more tangible if we restrict ourself to special covers by
imposing extra ramification conditions (e.g., by assuming that the order of 
the group is prime-to-$p$ or by imposing tameness conditions). 
In the case of $p$-group covers Raynaud has given a condition on the 
branch locus which 
eliminates vanishing cycles in the semi-stable models of the covers.

\begin{theorem}[\cite{Ra1}]\label{Raynaud} Let $X_K\to Y_K$ be a 
Galois cover 
with group $G$. Suppose 
$G$ is nilpotent and that $Y_K$ has a smooth model $Y$ such that the 
Zariski closure $B$ of the branch locus $B_K$ in $Y$ is \'etale over $R$ 
(i.e., the branch points
don't coalesce)
and that $X_K\to Y_K$ is tamely ramified (i.e., the inertia groups at points of 
$X_K$ are of order prime to char$(K)$). Then the special fiber of the stable 
model of $X_K$ 
is tree-like, i.e.\ the Jacobian of $X_K$ has potentially good reduction. 
 
\end{theorem}

The proof is an existence proof and it seems difficult to give an  explicit
one also in the simplest cases.  This is done in some sense for $p$-cyclic covers of 
the projective line in \cite{Le}, \cite{Ma}, \cite{Le-Ma4}.
There are still some questions 
relative to the monodromy and we note that little is known for 
$p^n$-cyclic covers ($n>1$) because 
Raynaud's condition on the branch locus doesn't transfer to intermediate 
covers. The complexity of the problem is partially measured by the wild 
monodromy.

\subsection{The wild monodromy}

Suppose $K$ is strictly henselian (with residue field $k$) and $X_K$ is a 
proper smooth curve over $K$ with $g(X_K)\geq 2$ or 
$X_K(K)\neq \emptyset $.
Then there exists an extension $L/K$ such that $X_L$ has  semi-stable reduction
and $L$ is contained in any extension $K'/K$ over which $X_{K'}$ has
semi-stable reduction.
Furthermore, the extension $L/K$ is Galois. 
If $g(X_K)\geq 2$ we denote by $\cal X$ 
the stable model of $X_L$. Then $\Gal (L/K)$ acts as a group of $K$-automorphisms on 
$X_L$  and this action extends as a group of automorphisms to $\cal X$ whose 
action on the special fiber ${\cal X}_k$ is faithful: 

$$\Gal (L/K) \hookrightarrow \Aut_k  {\cal X}_k \quad  (*)$$ 

With the above notation we use
the following terminology due to Raynaud (cf.\cite{Ra2} 2.2.2, and 4.2).

\begin{defi}\rm
The extension $L/K$ is called the {\it finite
monodromy extension}, its Galois group  $\Gal (L/K)$ is the 
{\it finite  monodromy} and its
 $p$-Sylow subgroup is the {\it wild monodromy}
 $\Gal (L/K)_w$.
\end{defi}


The following is a question asked by Colmez in the early 1990's.
\begin{qu}
Consider the hyperelliptic curve $Y^2=1-X^{2^n}$ over the $2$-adic
numbers. What is its stable reduction?
What is the wild monodromy and its action on the special fiber 
of the stable model? 
\end{qu}

\noindent 
Note that the branch locus of the hyperelliptic involution for $n>1$ doesn't satisfies Raynaud's 
condition. In fact one should regard the above curve as a $2^n$-cyclic
cover of the projective line ramified at the three points 
$Y=-1,1,\infty$. If we consider $Z=(Y+1)/2$ then we have the presentation
$X^{2^n}=1-Y^2=1-(2Z-1)^2=4(Z-Z^2)$. This is a $2^n$-cyclic
cover of the projective line ramified in the three points $0,1,\infty$ and
Raynaud's condition is satisfied.

\noindent 
\begin{qu}
Same question for $Y^4=(1+X)^2+X^3:=F(X)$.  This is a $2^2$-cyclic 
cover ramified in the 4 points $X=x_i, i=1,2,3$ and $X=\infty$ where 
$F(x_i)=0$; it satisfies Raynaud's condition.
\end{qu}

{\bf Now we concentrate on theorem \ref{Raynaud} and consider the case where 
the group is $G=\Z/p\Z$.}

\medskip

In \cite{Le-Ma4} we look at the case where the curve $X_K/G$ is $\PPP_K^1$ and we 
consider the following question: for a given genus $g=\frac{p-1}{2}(m-1)$ 
there are several types of degeneration and for each type there is an upper 
bound on the cardinality of the wild automorphism group 
$\Aut_k{\cal X}_k$ of the special fiber (see section 1). 

\begin{qu} \label{4}
For a given type 
of degeneration can the   
wild monodromy group $\Gal(L/K)_w$ be maximal (in the sense of attaining this bound)? 
\end{qu}

We have shown in \cite{Le-Ma4} that the answer is yes for
$p=2$ and $m=5$, i.e.\ genus 2 curves.

In other words, we are looking for a $p$-cyclic 
cover of the projective line over $\QQ_p^{\mbox{\rm \scriptsize tame}}$ which degenerates to the 
given type and for which 
the  wild monodromy group $\Gal(L/K)_w$ is a big as possible, i.e.\ attains
the upper bound. 
(It is not clear whether the group structure is always unique).
We have shown in \cite{Le-Ma4} that the answer to question \ref{4} is yes for
$p=2$ and $m=5$, i.e.\ genus 2 curves.

\noindent 
{In [Ma] there is an algorithm which produces a
polynomial whose zeroes are centers $y_i,i\in I$ of the blowing up which 
gives rise to the stable model. 
The wild monodromy extension is then essentially a sub\-extension of 
$K^{\mbox{\rm \scriptsize tame}}(y_i),i\in I$. The main 
problem with this algorithm is that the degree of the polynomial will be so 
big that we are 
not even able to write it down for $m\geq 9$. The size of $\Aut_k{\cal X}_k$ 
gives a bound 
for the possible size of the wild monodromy. 
We can use the results of \cite{Le-Ma1} (see section 1) together with the
injection (*)
in order to bound the degree and to say what the best algorithm is. 
This is what is done in \cite{Le-Ma4} and it is shown that we get the best 
algorithm at least in the case of potentially good reduction.

\subsection{Applications to the Inverse Galois Problem}

An affirmative answer to question \ref{4} would allow us to produce
finite monodromy extensions $L/K$ with big Galois groups, hence
Galois extensions of $K^{\mbox{\rm \scriptsize tame}}$ with big Galois groups. 

\begin{qu}
Is it possible to use this 
for the classical Inverse Galois problem in order to produce 
extensions with a given type 
of ramification at a given place?
\end{qu}

\noindent 
Let $K={\mathbb Q}_p^{\mbox{\rm \scriptsize tame}}$ and assume that
$C \la \PP_K$ is given birationally by the equation

$$Z^p=f(X)=1+cX^q+X^{q+1} \quad \mbox{\rm for }c\in R \mbox{ \rm and } q=p^s.$$

\begin{proposition}[\cite{Le-Ma4}]
a) If $v(c) \ge v(\lambda^{p/(q+1)})$ then $C$ has good reduction over $K$ 
and the wild monodromy is trivial.

\noindent b) If $c=p^{1/(q+1)}$ then $C$ has potentially good reduction and
the Galois group of the wild monodromy is  equal to the extraspecial 
group of order $pq^2$ and exponent $p$ if $p>2$ and 
equal to the central product $Q_8*D_8*...*D_8$ if $p=2$.
\end{proposition}
}

\begin{qu}
What is the ramification filtration
of the wild monodromy?
\end{qu}

\begin{qu}
How do the ramification filtrations of 
$\Gal (L/K)$ and $\Aut_k  {\cal X}_k$ behave under the injection (*)?
\end{qu}

\begin{qu}
A given type of degeneration corresponds to an analytic open in 
the moduli space $M_g$. The finite monodromy extension is locally constant 
(this is Krasner's Lemma).
What are the extensions of $K^{\mbox{\rm \scriptsize tame}}$ which can occur?
\end{qu}

{\small Claus Lehr,
Universit\`a di Padova,
Dipartimento di Matematica Pura ed Applicata,Via Belzoni 7, 35131 Padova, Italy, lehr@math.unipd.it

\medskip

 Michel Matignon,
Laboratoire de Th\'eorie des Nombres
et d'Algorithmique Arith\-m\'eti\-que,
UMR 5465 CNRS,
Universit\'e de Bordeaux I,
351 cours de la Lib\'eration, 
33405 Talence Cedex, France, matignon@math.u-bordeaux1.fr}

\chapter{Lifting Galois covers of smooth curves \\ {\mdseries \textsl{by} Michel Matignon}}

In this note we present some questions concerning the lifting 
of Galois covers of curves
from characteristic $p>0$ to characteristic zero.
We will focus on the case of elementary abelian $p$-groups which was 
studied by G. Pagot in his thesis.   

\section{Introduction}

Let $k$ be an algebraically closed field of characteristic $p>0$
and  $G$ a finite $p$-group.The group $G$ occurs as an automorphism group of 
$k[[z]]$  in many  ways; this is a consequence of the  Witt-Shafarevich theorem 
on the structure of the Galois group of a field $K$ of characteristic $p>0$.
This theorem asserts that the Galois group $I_p(K)$ of its maximal $p$-extension is  pro-$p$ 
free on $|K/\wp(K)|$ elements (as usual $\wp$ is the operator Frobenius minus identity). We apply this theorem to the power series field $K=k((t))$. 
Then $K/\wp(K)$ is infinite so we can 
realize $G$  in infinitely many ways as a quotient of $I_p$ and so as 
Galois group of a Galois extension $L/K$. The local field $L$ can be uniformized: namely 
the $t$-adic valuation of $K$ has a unique prolongation to a valuation 
$v_L$ of $L$ and if $z\in L $ is a uniformizing parameter ($v_L(z)=1$)
one has $L=k((z))$. If $\sigma \in G=\Gal (L/K)$, then $\sigma $ is an isometry of 
$(L,v_L)$ and so $G$ is a group of $k$-automorphisms of $k[[z]]$ with
fixed ring $k[[z]]^G=k[[t]]$.

When $G$ is abelian, one can use  Artin-Schreier-Witt theory in order to write down the extension $L/K$. 
For example in the case of  $G=\Z/p\Z$ the  $p$-cyclic extensions of $K$ are  
defined by the equation 

$(1)\quad w^p-w=f(t)\in K/\wp (K)$

\noindent
and $\sigma (w)=w+1$ is a generator for $G$.
In order to see $G$ as a $k$-automorphism group of the ring $k[[z]]$ we  need a desingularisation for $(1)$.
In this case it is easily obtained: from Hensel's lemma it follows that one can 
take $f(t)\in \frac{1}{t}k[\frac{1}{t}]$ with degree $m$ prime-to-$p$.
The integer $m$ is intrinsically defined and called the conductor of the cover; namely the $v_L$ valuation of the different ideal 
of the extension $L/K$ is $(p-1)(m+1)$. It is easy to see that then $w$ is 
an $m$th power in $L$. Further $z:=w^{-1/m}\in L$ is a uniformizing parameter
and $\sigma (z)=z(1+z^m)^{-1/m}$. It follows that there is up to change of parameter
one automorphism of order $p$ of $k[[z]]$ with conductor $m$ where $(m,p)=1$.

One can explicitly give a lifting of $\sigma$ as an automorphism 
of order $p$ of $\Z_p[\zeta][[Z]]$ where $\zeta$ is 
a primitive $p$-th root of unity, namely: $\Sigma (Z)=\zeta Z(1+Z^m)^{-1/m}$ works (the general case is easily deduced from $m=1$ case). There are other liftings as an automorphism group of  $R[[Z]]$ for 
some discrete valuation ring $R$  suitably ramified over $\Z_p[\zeta]$ which are not conjugate in the group $\Aut_RR[[Z]]$ to $\Sigma (Z)$. One way to look at the conjugacy class of an order $p$ 
automorphism of $R[[Z]]$, is to look at its set of fix points. The $p$-adic geometry of the open disc $\Spec R[[Z]]\otimes_R\Fr R$ induces a ``metric'' geometry on this set  and so we can attached to a conjugacy class 
a metric tree: the Hurwitz tree of the automorphism. 
The geometry of automorphisms of order $p$ of $R[[Z]]$ is now well understood through their 
Hurwitz tree (\cite{Gr-Ma2}, \cite{He} and \cite{Bo-We}) and so the problem of lifting $p$-cyclic actions to characteristic $0$ is satisfactorily understood. The analogous questions for a $G$-action when $G$ is not cyclic of order $p$ is far from being understood. 
In the sequel we try to make some questions. 
For a more general presentation,  the reader can have a look at  our MSRI notes (\cite{Ma2})  and at a recent survey by Q. Liu (\cite{Li}).

\section{ The local lifting problem} 

\begin{defi} The local lifting problem for a finite $p$-group action $G\subset \Aut_kk[\![z]\!]$
is to find a DVR,  $R$ finite over $W(k)$ and a commutative diagram

$$
\begin{array}{ccc}
{\rm Aut}_kk[\![z]\!] & \leftarrow&{\rm Aut}_RR[\![Z]\!]  \\
\uparrow &  \nearrow  & \\
G & & 
\end{array}
$$

A $p$-group $G$ has the local lifting property if the local lifting problem for all   actions
$G\subset \Aut_kk[[z]]$ has a positive answer.
\end{defi}

\noindent {In the introduction we have seen that a cyclic group of order $p$  has the local lifting property.  This is also the case for a cyclic group of order $p^2$  (\cite{Gr-Ma1}). However the proof
is more delicate because  the conjugacy classes of $p^2$-cyclic
actions  $G=\Z/p^2\Z \subset \Aut_kk[\![z]\!]$ are not parametrized by the conductors 
contrary to the $p$-cyclic actions. Another difficulty lies in the 
Sekiguchi-Suwa  deformation of the Artin-Schreier-Witt sequence 
to the Kummer sequence (\cite{Se-Su1}).
For $p^2$-cyclic \'etale covers it depends on tricky $p$-adic congruences and
we had to adapt their theory for each conjugacy class in order to produce 
an explicit lifting (over some DVR $R$) of the Artin-Schreier-Witt equations of the corresponding cover $k[[z]]/k[[z]]^G=k[[t]]$. This lifting is chosen in such a way that the generic fiber has a ``minimal singularity''. It is then shown that a desingularization of these relative 
$R$-curves produces a  (smooth) lifting of the action. This geometric method doesn't produce in an 
explicit way an automorphism of order $p^2$ of $R[[Z]]$.  

For $p^n$-cyclic actions ($n>2$) T. Sekiguchi and N. Suwa 
have a satisfactory 
generalization of their theory  (see \cite{Se-Su2}, \cite{Se-Su3}), so it is reasonable to expect a positive answer of 
the local lifting problem for all $p^n$-cyclic actions (also called Oort's conjecture). 

Recently, G. Pagot (\cite{Pa2}, \cite{Pa3}, \cite{Ma3}) has shown that the group $(\Z/2\Z)^2$ has the local lifting property ($p=2$). 

In general, for a non cyclic action,  
there are combinatorial obstructions (see \cite{Be}, \cite{Gr-Ma1}) or obstructions
of differential flavour (\cite{Pa1}, \cite{Pa2}, see section 4
below).}

\section{ Inverse Galois local lifting problem for $p$-groups.} 

Let $G$ be a finite $p$-group; in the introduction we saw that $G$ occurs as 
a group of
$k$-automorphism 
of $k[[z]]$ in many ways, so we would like to address a weaker problem than the local 
lifting problem.

\begin{defi} For a finite $p$-group $G$ we say that $G$ has the weak local lifting property if 
there exists a DVR, $R$ finite over $W(k)$, a faithful action $i: G\rightarrow \Aut_kk[[z]]$ and a 
commutative diagram

$$
\begin{array}{cccc}
{\rm Aut}_kk[\![z]\!] & \leftarrow&{\rm Aut}_RR[\![Z]\!]  \\
i\uparrow &  \nearrow  & \\
G & & 
\end{array}
$$
\end{defi}

In \cite{Gr-Ma2} we prove that $p^n$-cyclic groups have the weak local lifting property.

In \cite{Ma1} we prove that elementary abelian $p$-groups have the weak local lifting property.

\begin{qu} Let $G$ be a non abelian group of order $p^3$. Does $G$ have
the weak local lifting property?
\end{qu}

\section{ Local lifting problem for elementary abelian $p$-groups after Pagot's thesis}

In his thesis G. Pagot studies actions of elementary abelian $p$-groups on
$k[[z]]$ which can be lifted to automorphism groups of $R[[Z]]$ for some 
DVR $R$ finite over $W(k)$. In this context there remain various
natural questions to solve. We would like to mention some.

\begin{defi}\rm  We fix an infinite point $\infty \in \PPP^1_k$. A space $L_{m+1,n}$
is a $\F_p$-vector space of dimension $n$ whose non-zero elements are logarithmic
differential forms on $ \PPP^1_k$ with only one zero which is at $\infty$ and 
of order $m-1$. 
\end{defi}

Necessarily $(m,p)=1$. Such a space naturally occurs when considering 
$(\Z/p\Z)^n$ actions on $R[[Z]]$ and conversely to such a space one can 
naturally associate  $(\Z/p\Z)^n$ actions on $R[[Z]]$ (\cite{Ma1}, \cite{Pa1}).

\begin{qu}
  Give necessary and sufficient conditions on $m,n$ for the existence of spaces $L_{m+1,n}$. In particular prove or disprove the following:
assume that $n\geq 2$ and that there is a space $L_{m+1,n}$; then $p^{n-1}(p-1)$ divides $m+1$. 
\end{qu}

\medskip
\noindent 
{\bf Case  $n=1$.} 

\medskip

\noindent
Let $m\in \N-p\N$; we can describe the spaces $L_{m+1,1}$ (see [Gr-Ma2], [He]): These   are 
the spaces $\F_p df/f$, where $f=\prod_{1\leq i\leq m+1}(z-x_i)^{h_i}$
satisfies the following conditions:
\begin{enumerate}
\item $\sum_{1\leq i\leq m+1}h_ix_i^{\el}=0 $ for $1\leq \el \leq m-1$

\item $\prod_{i<j}(x_i-x_j)\neq 0$

\item $x_i\in k,h_i\in \Z-p\Z$.
\end{enumerate}

\begin{example} Let $m\notin pZ$ and $f=z^{-m}-1$. Then  $$\frac{df}{f}= \frac{mdz}{z^{m+1}-z}$$ generates an 
 $L_{m+1,1}$.
\end{example}

\begin{example} Let $m=p-1$ and  $f=\prod_{1\leq i\leq p-1}(z-i)^i$. Then $$\frac{df}{f}=\frac{dz}{1-z^{p-1}}$$ generates an $L_{p,1}.$
\end{example}

\noindent 
{\bf Case $n=2$.}

\medskip

\noindent
G. Pagot gives the following characterization.

\begin{proposition}[\cite{Pa1}] Let $\omega_0,\omega_p$ be two differential forms over $\PPP^1_k$.
Then $\F_p \omega_0+\F_p \omega_p$ is an $L_{m+1,2}$
iff $p|m+1$  and there are 2 polynomials $A$ and $B$ with
$$\forall [i,j] \in \PPP^1 (\F_p)\quad \deg (iA+jB)=(m+1)/p,$$
$$((A^p-AB^{p-1})^{p-1})^{(p-1)}=-1$$ 
and 
$$\omega_0=\frac{Adz}{A^pB-AB^p}, \ \omega_p=\frac{Bdz}{A^pB-AB^p}.$$
\end{proposition}

Finding  $A$ and $B$ as above is a difficult problem. In the small degree case G. Pagot was able to deduce the following theorem which justifies the question above:

\begin{theorem}[\cite{Pa1}] Let $p>2$. 

-There is no space $L_{p,2}$;

-If there is a space $L_{2p,2}$, then $p=3$;

-There is no space $L_{3p,2}$.

\end{theorem}

Note that G. Pagot has described the spaces $L_{6,2}$ when $p=3$ and the spaces $L_{m+1,2}$ 
when $p=2$.

\medskip
\noindent
{\bf Case $n>2$.}

\medskip

\noindent{
Let  $(\omega_1,...,\omega_n)$ be a basis for an  $L_{m+1,n}$. Then $p^{n-1}|m+1$ 
and these $n$ forms have  $(m+1)(p-1)^{n-1}/p^{n-1}$ poles in common (\cite{Pa1}).
Let $m+1=p^{n-1}(p-1)$ and
$a_1,...,a_n\in \F _p^{alg}$, pairwise distinct.
For $1\leq j\leq n$, let 
 $$f_j=\prod_{(\epsilon_1,...,\epsilon_n)\in \{0,1,...,p-1\}^n}(z-\sum_{1\leq i\leq n}\epsilon_ia_i)^{\epsilon_j}$$ and  $\omega_j=df_j/f_j$. Then 
there exists $P\in \F_p[x_1,...,x_n]$ such that if $$P(a_1,...,a_n)\neq 0$$ 
then $\sum_{1\leq j\leq n}\F_p\omega_j$ is an $L_{m+1,n}$ (see \cite{Ma1} and \cite{Pa1}). }

\begin{qu} 
Assume $p=2$. Consider an action of  $G=(\Z/2\Z)^n$ as an automorphism group of $k[[z]]$ which satisfies Bertin's numerical conditions (see 
\textup{\cite{Be}, \cite{Ma3}}). Is it possible to lift such an action as an automorphism 
group of $R[[Z]]$ for some DVR finite over $\Z_2$? 
\end{qu}

\noindent
 If $n= 2$, the numerical conditions are empty; G.  Pagot
shows  that the answer is yes (\cite{Pa2}, \cite{Pa3}). In a handwritten paper (2004) he answers positively 
to the case $n=3$ and it seems he can show that his proofs generalize to 
the general case $n>3$.

\bigskip

{\small Michel Matignon,
Laboratoire de Th\'eorie des Nombres
et d'Algorithmique Arith\-m\'eti\-que,
UMR 5465 CNRS,
Universit\'e de Bordeaux I,
351 cours de la Lib\'eration, 
33405 Talence Cedex, France, matignon@math.u-bordeaux1.fr}

\chapter{Abelian varieties isogenous to a Jacobian \\ {\mdseries \textsl{by} Frans Oort}}

\makeatletter
\renewcommand{\subsection}[1]{\@startsection{subsection}{2}{\z@}
            {-3.25ex plus -1ex minus -.2ex}{-1sp}{\normalsize\bf}
             {\ignorespaces#1 }}
\renewcommand{\thesubsection}{(\thesection.\arabic{subsection})} 
\makeatother

\section*{Introduction}
\subsection{} {\bf Question.}  {\it Given an abelian variety $A$; does there exist an algebraic curve $C$ such that there is an isogeny between $A$ and the Jacobian of $C$} ?
\begin{itemize}
\item If the dimension of $A$ is at most three, such a curve exists; see \ref{<4}.
\item For any $g \geq 4$ there exists an abelian variety $A$ of   $\dm(A)= g$ over $\CC$ such that there is no  algebraic curve $C$ which admits an isogeny $A \sim \Jac(A)$, see \ref{CC}. One of the arguments which proves this fact (uncountability of the ground field) does not hold over a countable field. 
\end{itemize}
\n
Therefore:
\begin{itemize}
\item The question remains open  over a countable  field, see \ref{E1}. One can expect that the answer to the question in general is negative for abelian varieties of dimension $g \geq 4$ over a given field.
\item We offer a possible approach to this question via Newton polygons in positive characteristic, see \ref{E3}.
\end{itemize}

\section{Jacobians and the Torelli locus}

\subsection{\bf Jacobians.} Let $C$ be a complete curve over a field $K$. We write $J(C) = \Pic^0_{C/K}$.

\vn
In case $C$ is irreducible and non-singular we know that $J(C)$ is an abelian variety.  Moreover, $J(C)$ has a canonical polarization. This {\it principally polarized abelian variety} $\Jac(C) = (J(C), \Theta_C = \lambda)$ is called the {\it Jacobian} of $C$. 

\vn 
Suppose $C$ is a geometrically connected, complete curve of genus at least 2 over a field $K$. We say that $J(C)$ is a  {\it curve of compact type}  if $C$ is  a  stable curve such that:

-  its geometrically irreducible components are non-singular, and

-  its dual graph has homology equal to zero;

\vn
equivalently: $C$ is stable, and for an algebraic closure $k \supset K$ the curve $C_k$ is a tree of non-singular irreducible components;

\vn
equivalently (still $g \geq 2$): $C$ is a stable curve, and $J(C)$ is an abelian variety.

\vn
For $g=1$ we define ``of compact type''  as ``irreducible + non-singular''.

\vn
The terminology ``a curve of compact type'' is the same as ``a good curve'' (Mumford), or ``a nice curve''. 

\subsection{\bf The Torelli locus.} By $C \mapsto \Jac(C) := (J(C), \Theta_C = \lambda)$ we obtain a morphism $j: \cM_g \to \cA_{g,1}$, from the moduli space of curves of genus $g$ to the moduli space of principally polarized abelian varieties; this is called {\it the Torelli morphism}. The image 
$$\cM_g \twoheadrightarrow \cT^0_g \to \cA_{g,1}$$ 
is called the {\it open Torelli locus}.

\vn
Let  $\cM^{\sim}_g$ be the moduli space of  curves of compact type. The Torelli morphism can be extended to a morphism $\cM^{\sim}_g \to \cA_{g,1}$;  its image  

$$\cM^{\sim}_g \twoheadrightarrow \cT_g \to \cA_{g,1}, \quad\quad \cT_g = (\cT_g^0)^{\rm Zar},$$
is the Zariski closure of $\cT^0_g$; we say that  $\cT_g$ is  the {\it closed Torelli locus}.

\vn
From now on let $k$ be an algebraically closed field.

\subsection{}\label{<4} {\it Suppose $1 \leq g \leq 3$. Then  every abelian variety $A$ of dimension $g$ is isogenous with the Jacobian of a curve of compact type.} 
 \\
{\bf Proof.} In fact,  a polarized abelian variety $(A,\lambda)$ over an algebraically closed field is isogenous with a principally polarized abelian variety $(B,\mu)$. We know that there exists a  curve $C$ of compact type with $(B,\mu) \cong \Jac(C)$; for $g=1$ this is clear; for 
$g=2$ see A. Weil, \cite{Weil}, Satz 2; for $g=3$ see F. Oort \& K. Ueno, \cite{FO.U}, Theorem 4.
\B

\subsection{} Let $g \in  \ZZ_{>0}$. We write $Y^{(cu)}(k,g)$ for the statement:

\vn
$Y^{(cu)}(k,g)$. {\it There exists an abelian variety $A$ defined over $k$ such that there does not exist a  curve $C$ of compact type of genus $g$ defined over $k$ and an isogeny $A \sim J(C)$} (here $c$ stands for ``of compact type'', and $u$ stands for ``unpolarized'').

\section{Other formulations}
\subsection{}
$Y^{(cp)}(k,g)$. {\it There exists a polarized  abelian variety $(A,\lambda)$ with $\dim(A) = g$ defined over $k$ such that there does not exists a  curve $C$ of compact type (of genus $g$) defined over $k$ and an isogeny $(A,\lambda) \sim \Jac(C)$.}

\subsection{}
$Y^{(iu)}(k,g)$. {\it There exists an abelian variety $A$  with $\dim(A) = g$ defined over $k$ such that there does not exists an irreducible curve $C$ (of genus $g$) defined over $k$ and an isogeny $A \sim J(C)$.}

\subsection{}
$Y^{(ip)}(k,g)$.  {\it There exists a polarized  abelian variety $(A,\lambda)$  with $\dim(A) = g$ defined over $k$ such that there does not exists an irreducible curve $C$ (of genus $g$) defined over $k$ and an isogeny $(A,\lambda) \sim \Jac(C)$.}

\subsection{} Note that $Y^{(cu)} \Rightarrow Y^{(cp)} \Rightarrow Y^{(ip)}$ and $Y^{(cu)} \Rightarrow Y^{(iu)} \Rightarrow Y^{(ip)}$.

\vn
Given a point $[(A,\lambda)] = x \in \cA_g$ in the moduli space of polarized abelian varieties we write $\cH(x) \subset \cA_g$ for the {\it Hecke orbit} of $x$; by definition $[(Y,\mu)] = y \in \cH(x)$ if there exists an isogeny $A \sim B$ which maps $\lambda$ to a rational multiple of $\mu$.

\subsection{} Here is a reformulation:
$$Y^{(ip)}(k,g) \ \Longleftrightarrow\ \cH(\cT^0_g)(k) \subsetneqq \cA_g(k),\
\mbox{\rm where} \ \cH(\cT^0_g) = \cup_{x \in \cT^0_g} \ \  \cH(x).$$

\section{Over large fields}

\subsection{}\label{CC} {\it Suppose $g \geq 4$ and let $k$ be an algebraically closed field which is uncountable, or a field such that ${\rm tr.deg.}_P(k) > 3g-3$ (here $P$ is the prime field of $k$). Then $Y^{(cu)}(k,g)$ holds.} For example if $k = \CC$ we know that there is an abelian variety of dimension $g$ not isogenous to a Jacobian variety of any curve of compact type. 

\subsection{}\label{CC2}
{\it We show that $Y^{(cp)}(k,g)$ holds for $k=\CC$ and $g \geq 4$}.\\
{\bf Proof.} In this case $\dm(\cM_g \otimes \CC) = 3g-3 < g(g+1)/2 = \dm(\cA_g \otimes k)$. Hence $\cT_g \otimes k$ is a proper subvariety of $\cA_g \otimes k$. Write $\cH(\cT_g \otimes k)$ for the set of points corresponding with all polarized abelian varieties isogenous with a (polarized) Jacobian ($\cH$ stands for  ``Hecke orbit''). We know that $\cH(\cT_g \otimes k)$ is a countable union of lower dimensional subvarieties. Hence $\cH(\cT_g(\CC)) \subsetneqq \cA_g(\CC)$. 
\B 

\vn
An analogous fact can be proved in positive characteristic using the fact that Hecke correspondences are finite-to-finite on the ordinary locus, and that the non-ordinary locus is closed and has codimension  one everywhere.

\subsection{} A referee asked whether I could give an example illustrating \ref{CC2}. Let $V \subset \cA_g \otimes \overline{\QQ}$ be an irreducible subvariety with $\dim(V) > 3g-3$; for example choose $V$ to be equal to an irreducible component of $\cA_g \otimes \overline{\QQ}$. Let $\eta$ be its generic point, $K := \overline{\QQ}(\eta)$, with algebraic closure $\overline{K} = L$; choose an embedding $L \hookrightarrow \CC$. Over $L$ we have a polarized abelian variety $(A,\lambda)$ corresponding with $\eta \in 
\cA_g(L)$. This abelian variety is not isogenous with the Jacobian of a curve.  However, I do not know a ``more explicit example''.

\vn 
The previous proof uses the fact that the transcendence degree of $k$ is large. For a field like $\overline{\QQ}$ this proof cannot be used. However we expect the following to be true.

\subsection{}\label{E1}
{\bf Expectation} (N. Katz). {\it One can expect  that  $Y^{(cu)}(\overline{\QQ},g)$ holds for every $g \geq 4$.}  For a more general question see \cite{Milne}, 10.5.\\
{\it  One can expect  that  $Y^{(cu)}(\overline{\FF_p},g)$ holds for every $g \geq 4$ and every prime number $p$.} 

\vn
{\bf Remark.}  
$Y^{(cu)}(\overline{\FF_p},g) \ \Rightarrow Y^{(cu)}(\overline{\QQ},g); \quad\mbox{see the proof of \ref{p0}}.$

\section{Newton polygons and the $p$-rank}
\n
In order to present an approach to \ref{E1} we recall some notions.

\subsection{} Manin and Dieudonn\'e proved that isogeny classes of $p$-divisible groups over an algebraically closed field are classified by their Newton polygons, see \cite{Manin}, page 35. A symmetric Newton polygon (for height $h = 2g$) is a polygon in $\QQ \times \QQ$:  

- starting at $(0,0)$, ending at $(2g,g)$, 

- lower convex,

- having breakpoints in $\ZZ \times \ZZ$, and

- a slope $\lambda \in \QQ$, \ \ $0 \leq \lambda \leq 1$, appears with the same multiplicity as $1-\lambda$.\\
See \cite{FO-EO}, 15.5. 

\vn
An abelian variety $A$ in positive characteristic determines a Newton polygon $\cN(A)$ by taking the Newton polgon of $X = A[p^{\infty}]$. This defines a symmetric Newton polygon.

\vn 
The finite set of symmetric Newton polygons belonging to $h=2g$ is partially ordered by saying that $\xi' \prec \xi$ if no point of  $\xi'$ is below  $\xi$, colloquially: if $\xi'$ is ``above'' $\xi$. 

\subsection{} Let $\xi$ be a symmetric Newton polygon. We write 
\begin{eqnarray*} & & W_{\xi} = \{[(B,\mu)]\mid \cN(B) \prec \xi\} \subset\cA_{g,1}; 
\\ & & W^0_{\xi} = \{[(B,\mu)]\mid \cN(B) = \xi\} \subset\cA_{g,1}.
\end{eqnarray*}
Grothendieck proved that under specialization Newton polygons go up
and  Grothendieck and Katz showed that the locus $W_{\xi} \subset \cA_{g,1} \otimes \FF_p$ is {\it closed}, see \cite{Katz}.  The locus $W^0_{\xi} \subset \cA_{g,1} \otimes \FF_p$ is {\it locally closed} and  $\overline{W^0_{\xi}} = W_{\xi}$.

One can also define  $W_{\xi}$ for $\cA_g \otimes \FF_p$;  these loci are  closed; we will focus on the principally polarized case.

\vn
These loci $W_{\xi}$ are now reasonably well understood in the principally polarized case. The codimension of $W_{\xi}$ in $\cA : = \cA_{g,1} \otimes \FF_p$ is precisely the length of the longest chain from $\xi$ to the lowest Newton polygon $\rho$ ($g$ slopes equal to 0, and $g$ slopes equal to 1: the ``ordinary case''), see \cite{FO-CH} and \cite{FO-NP}. In particular the highest Newton polygon $\sigma$ (all slopes equal to 1/2: the  ``supersingular case'') corresponds to a closed subset of dimension $[g^2/4]$ (conjectured by T.Oda \& F. Oort; proved by K.-Z. Li \& F. Oort, \cite{Li.FO}, and reproved in \cite{FO-NP}). In particular:  {\it the longest chain of symmetric Newton polygons is equal to} 
$g(g+1)/2 - [g^2/4]$.

\subsection{} For an abelian variety $A$ over an algebraically closed field $k \supset \FF_p$ we define the $p$-rank $f(A)$ of $A$ by: 
$$A(k)[p] \quad\cong\quad (\ZZ/p)^{f(A)}.$$
Here $G[p]$ for an abelian group $G$ denotes the group of $p$-torsion points. It is easy to see that all values $0 \leq f \leq g$ do appear on $\cA_g \otimes \FF_p$.

\subsection{\bf Intersection with the Torelli locus.} We study intersection  $W_{\xi} \cap (\cT^0_g \otimes \FF_p)$ and the intersection $W_{\xi} \cap (\cT_g \otimes \FF_p)$. In low dimensional cases, and in some particular cases the dimension of these intersections is well-understood. However,
in general these intersections are difficult to study. Some explicit cases show that in general the dimension of an irreducible component of  $W_{\xi} \cap (\cT^0_g \otimes \FF_p)$ need not be equal to    $\dm(W_{\xi}) + \dm(\cM_g) - \dm(\cA)$.

\subsection{}{\bf The $p$-rank.}
We write 
$$V_f = \{[(A,\lambda)] \in \cA_g \otimes \FF_p \mid f(A) \leq f\}.$$
This is called a $p$-rank stratum.  We know:
$$\dm(V_f) \quad=\quad g(g+1)/2 - (g-f).$$
For  principally polarized abelian varieties this was proved by Koblitz; the general case can be found in \cite{N.FO}.

\subsection{}{\bf Remark.} Every $0 \leq f \leq g$ there exists a Newton polygon $\xi$ such that $V_f = W_{\xi}$. In other words, the Newton polygon stratification refines the $p$-rank stratification.

We see in \cite{F.vdG} that the dimension of every component of $V_f \cap (\cM_g \otimes \FF_p)$ equals $3g-3-(g-f)$ (for $g>1$).

\section{From positive characteristic to characteristic zero} 
In this section we formulate a  question; a positive answer to this would imply that $Y^{(\cdot \cdot)}$ holds over $\overline{\QQ}$. Here is the argument showing this last statement:

\subsection{}\label{p0} {\it Suppose given $g$ and a Newton polygon $\xi$} 
(of height $2g$).  
\begin{center}
$W^0_{\xi} \cap (\cT_g \otimes \FF_p) = \emptyset$ \quad$\Rightarrow$\quad $Y^{(cu)}(\overline{\QQ},g)$.
\end{center} 
{\bf Proof.} The condition $W^0_{\xi} \cap (\cT_g \otimes \FF_p) = \emptyset$  means that $\xi$ does not appear on $\cM^{\sim}_g$; hence there is an abelian variety $A_0$ with  $\cN(A_0) = \xi$ over $\overline{\FF_p}$,   which is not isogenous with the Jacobian of a  curve of compact type over $\overline{\FF_p}$. Choose an abelian variety $A$ over $\overline{\QQ}$ which has good reduction at $p$, and whose reduction is isomorphic with $A_0$ (this is possible by P. Norman \& F. Oort, see \cite{N.FO}). We claim that the abelian variety $A$ satisfies the condition  $Y^{(cu)}(\overline{\QQ},g)$: a  curve $C$ of compact type  over $\overline{\QQ}$ with $A \sim_{\overline{\QQ}} J(C)$ would have a Jacobian  $J(C)$ with good reduction $J(C)_0 \sim A_0$; this shows that $C$ has compact type reduction, and that $J(C_0) = J(C)_0 \sim_{\overline{\FF_p}} A_0$; this is a contradiction; this proves the implication.
\B

\subsection{} Note that for $g>1$ we have:
$$g(g+1)/2 - [g^2/4] < 3g-3 \quad\Longleftrightarrow\quad g \leq 8;$$  $$g(g+1)/2 - [g^2/4] > 3g-3 \quad\Longleftrightarrow\quad g \geq 9.$$

\subsection{\bf Expectation.}\label{Ep}  {\it Let $g = 11$, and let $\xi$ be the Newton polygon with slopes $5/11$ and $6/11$. We 
expect:
$$W^0_{\xi} \cap (\cT_{11} \otimes \FF_p) \quad\stackrel{?}{=}\quad \emptyset,$$
i.e.\ we think that this Newton polygon should not appear on the moduli space of  curves of compact type of genus equal to 11.} 

\vn
More generally one could consider $g \gg 0$, and $\xi$ given by slopes $i/g$ and $(g-i)/g$ such that gcd$(i,g)=1$, and such that the codimension of $W_{\xi}$ in $\cA_g$ is larger than $3g-3$:

\subsection{\bf Expectation.}\label{E3}  {\it  Suppose given $g$ and a Newton polygon $\xi$ 
(of height $2g$). Suppose:}
\begin{itemize}
\item {\it the longest chain connecting $\xi$ with $\rho$ is larger than $3g-3$}; 
\item
{\it the Newton polygon has ``large denominators''.}
\end{itemize} 
{\it Then we expect that}
 $$W^0_{\xi} \cap (\cT_g \otimes \FF_p) \stackrel{?}{=} \emptyset.$$
Note that the first condition implies ($g=8$ and $\xi = \sigma$) or $g>8$.
We say that the Newton polygon has {\it large denominators} if all slopes written as rational numbers with coprime nominator and denominator have a large denominator, for example at least eleven.

If a Newton polygon $\xi$ does not appear on $\cT_g$ one could expect that the set of slopes of $\xi$ is not a subset of slopes of a Newton polygon appearing on any $\cT_h$ with $h \geq g$.

\subsection{} We do not have a complete list for which values of $g$ and of $\xi$ we have $W^0_{\xi} \cap \cT_{g} \not= \emptyset$, not even  for relatively small values of $g$.

One could also study which Newton polygons show up on the open Torelli locus $\cT^0_g \otimes \FF_p$.

\subsection{} Probably there is a genus $g$ and a symmetric Newton polygon $\xi$ for that genus such that $\xi$ does show up on $\cT_g \otimes \FF_p$  and such that  $Y^{(nu)}(\overline{\QQ},g)$ is true. In other words: it might be that our proposed attempt via \ref{E3} can confirm \ref{E1} for large $g$, but not for all $g \geq 4$. 

\subsection{}{\bf Conjecture.}\label{E4}  {\it Let $g', g'' \in \ZZ_{>0}$; let $\xi'$, respectively $\xi''$ be a symmetric Newton polygon appearing on $\cT^0_{g'} \otimes \FF_p$, respectively on  $\cT^0_{g''} \otimes \FF_p$; write $g=g'+g''$. Let $\xi$ be the Newton polygon obtained by taking all slopes with their multiplicities appearing in $\xi'$ and in $\xi''$. We conjecture that in this case $\xi$ appears on $\cT^0_{g}$.}

\subsection{} We give some references. 

\vn
In \cite{vdG.vdV1} the authors show that for every $g \in \ZZ_{>0}$ there exist an irreducible curve of genus $g$ in characterstic 2 which is supersingular.  One can expect that for every positive $g$ and every prime number $p$ there exists an irreducible curve of genus $g$ in  characterstic $p$ which is supersingular; this would follow if \ref{E4} is true.  For quite a number of values of $g$ and $p$ existence of a supersingular curve has been verified, see \cite{Re}, Th. 5.1.1.

\vn
Next one can ask which Newton polygons show up on the hyperelliptic locus $H_g$. Here is a case where that dimension is known: {\it for $g=3$ every component of $W_{\sigma} \cap (H_3 \otimes \FF_p)$ has the expected dimension $5 + 2 - 6 = 1$}, see \cite{FO-H}; however already in this ``easy case'' the proof is quite non-trivial.

\vn
In \cite{SZ1} the authors show that for a hyperelliptic curve in characteristic two of genus $g = 2^n-1$ and $2$-rank equal to zero, the smallest slope equals $1/(n+1)$. In \cite{SZ2} we see which smallest slopes are possible on the intersection of the hyperelliptic  locus with $V_{g,0}$  for $g<10$.

\vn
In \cite{G.P1} it is shown that components of the intersection of $V_{g,f}$ with the hyperelliptic locus all have dimension equal to $g-1+f$ (i.e. codimension $g-f$ in the hyperelliptic locus). In \cite{G.P2} we find the question whether this intersection is transversal at every point.

\vn
Instead of the Newton polygon stratification one can consider another stratification, such as the the ``Ekedahl-Oort stratification'' or the ``stratification by $a$-number''. In  \cite{Re}, Chapter 2 (especially Th. 2.4.1, also see the remark at the end of that chapter) we see that for an $a$-number $m$ and a prime number $p$ for every $g$ with $g > pm + (m+1)p(p-1)/2$ the $a$-number $m$ does not appear on $\cT^0_g \otimes \FF_p$. As the  $a$-number is not an isogeny invariant, this fact does not contribute directly to the validity of a statement like \ref{E3}; however it does show that intersecting a stratification on $\cA_g \otimes \FF_p$ with the Torelli locus presents difficult and interesting problems.

\bigskip

{\small Universiteit Utrecht, Mathematisch Instituut, Postbus 80.010, 3508 TA Utrecht, Nederland, oort@math.uu.nl}

\chapter{Minimal maximal number of automorphisms of curves \\ {\mdseries \textsl{by} Frans Oort}}

\begin{defi} Let $k$ be an algebraically closed field, and let $g \in \ZZ_{\geq 2}$. We write:
$$\mu(g,{\rm char}(k)) \quad:=\quad {\rm Max}_C \ \ \#(\Aut(C)),$$
the maximum taken over all complete, irreducible and nonsingular algebraic curves $C$ of genus $g$ defined over $k$.
\end{defi}

Note that once a characteristic is chosen, the number $\mu(g,{\rm char}(k))$ does not depend on the choice of $k$.

\begin{thm}[Characteristic zero] 
The following is know if ${\rm char}(k)=0$: 
\begin{enumerate} 
\item (Macbeath \textup{\cite{Mac}}) The Hurwitz upper bound $\mu(g,0) \leq   84{\cdot}(g-1)$ is attained for infinitely many values of $g$.
\item (Accola \textup{\cite{Acc}, \cite{Acc2}} and Maclachlan \textup{\cite{Macl}}) 
\begin{itemize}
\item  $\mu(g,0) \geq 8(g+1)$, for every $g \geq 2$;
\item  for infinitely many values of $g$ we have $\mu(g,0) = 8(g+1)$. 
\end{itemize}
\end{enumerate}
\end{thm}

\begin{question} Suppose given a prime number $p$. Does there exist a polynomial $M_p \in \QQ[T]$ such that:
\begin{itemize}
\item  $\mu(g,p) \geq M_p(g)$, for every $g \geq 2$, and 
\item  for infinitely many values of $g$ we have $\mu(g,p) = M_p(g)$ ? \end{itemize}
\end{question}

\begin{remark}
If a polynomial with these properties exists, it is unique.
\end{remark}

\begin{question}
If a polynomial such as asked for in the previous question exists, is it the same for different prime numbers? Is it equal to $8(g+1)$?
\end{question}

\bigskip

{\small Universiteit Utrecht, Mathematisch Instituut, Postbus 80.010, 3508 TA Utrecht, Nederland, oort@math.uu.nl}

\end{document}